\documentclass[a4paper,12pt]{article}
\usepackage{amsmath,amsthm}
\usepackage{txfonts}
\usepackage{enumerate}

\numberwithin{equation}{section}

\newtheorem{thm}{Theorem}[section]

\theoremstyle{definition}

\theoremstyle{remark}

\def\P{{\mathsf{P}}}

\def\R{ {\mathbf R} }

\def\Z{ {\mathbf Z} }

\def\Ex{ \mathbf{E} }

\def\P{ \mathsf{P} }

\title{\bf Discrete It\^o Formulas and Their Applications to 
Stochastic Numerics}
\author{Jir\^o Akahori \\[12pt]
Department of Mathematical Sciences \\
Ritsumeikan University \\[8pt]
1-1-1, Nojihigashi, Kusatsu, Shiga, 525-8577, Japan \\[6pt]
E-mail: akahori@se.ritsumei.ac.jp}
\date{}

\begin{document}
\maketitle





\section{Introduction}
This is a survey note of the author's observations on 
the discrete-time analogues of  It\^o formulas.
The observations are summarized 
as follows. 
\begin{enumerate}
\item The standard It\^o formula
is based on the stochastic integrals, 
while in discrete-time frameworks, 
we instead rely on (conditional) Fourier series expansions. 
Detailed explanations will be given in Section \ref{DIF}. 

\item In a parallel way 
that the standard one describes the Kolmogorov equation 
for a given stochastic differential equation, 
the discrete It\^o formula gives a finite difference equation 
for a given approximating equation (Euler-Maruyama, for example) 
of SDE. 
This observation leads to a computational 
framework of Monte-Carlo simulations of the 
finite difference scheme for partial differential equations
in a high dimension.
(Section \ref{high})


\end{enumerate}

\section{Discrete It\^o formulas}\label{DIF}
\subsection{Fujita's It\^o formula}
Let us start with Fujita's Discrete It\^o Formula (DIF for short) 
\cite{Fujita} 
for the simple random walk;
\begin{equation}\label{simple}
W_k = \tau_1 + \tau_2 + \cdots + \tau_{k},
\end{equation}
where $ \{ \tau_1,...,\tau_t,...\} $ is a Bernoulli sequence such that
$ \P (\tau_k =  \pm 1 ) = 1/2 $.
By a simple algebra, we have
\begin{equation}\label{FDIF}
F( \tau_k ) = \frac{\{ F (+1) - F(-1) \}}{2} \tau_k + 
\frac{\{ F (+1) + F(-1) \}}{2}
\end{equation}
for every $ F : \{ -1, +1 \} \to \R $.
A DIF for $ f : \Z \to \R $ is obtained 
by regarding $ f (W_{t+1}) = f (W_t + \tau_t) $ as
a function on $ \{ -1, +1 \} $ and applying (\ref{FDIF}) to it:
\begin{equation*}
\begin{split}
& f ( W_{t+1} ) - f (X_t) = f (W_t + \tau_{t+1}) - f (W_t) \\
&= \frac{f ( W_t + 1) -f( W_t - 1) }{2} 
\tau_{t+1} + \frac{f ( W_t + 1) +f ( W_t - 1) }{2} - f( W_t).
\end{split}
\end{equation*}
The starting point is to regard (\ref{FDIF}) as
Fourier expansion of $ F $ with respect to the orthonormal basis 
$ \{ 1, \tau_k \} $.

\subsection{The first generalization} 
Let $ \xi $ be a real valued random variable with $ \Ex [\xi] = 0 $ 
and $ \mathrm{Var} [\xi] = \Ex[\xi^2 ] < \infty $.
Let $ \nu $ be its law and choose
an orthonormal basis $ \{ H_n : n \in \mathbf{N} \} $
of the Hilbert space $ L^2 (\R, \nu) $
by expanding $ H_0 \equiv 1, H_1 (x) = x/\sqrt{\mathrm{Var} [\xi]} $.
Then the DIF for the random walk (the sum of independent copies of $ \xi $)
\begin{equation*}
W_{t_k} = \xi_1 +\xi_2 +\cdots +\xi_k, \quad (k \in \Z)
\end{equation*}
would be the following orthogonal expansion of 
$ f (W_{t_{k+1}} + \,\,\cdot \,\,) 
- f (W_{t_k} ) $. 
\begin{equation}\label{ADIF1}
\begin{split}
f ( W_{t_{k+1}} ) - f (W_{t_k}) 
&= f (W_{t_k} + \xi_{k+1}) - f (W_{t_k}) \\
&= \frac{1}{\Ex [|\xi|^2 ]}
\left( \int  f ( W_{t_k} + x) x \,\nu(dx) \right) 
(W_{t_{k+1}} - W_{t_k})  \\
&\quad + \frac{1}{t_{k+1} -t_k}
\left( 
\int \{ f(W_{t_k} +x) - f(W_{t_k})\} \,\nu(dx)
\right)
\,(t_{k+1} -t_k) \\
&\qquad +\sum_{n=2}^\infty \left( \int f ( W_{t_{k+1}} + x) 
H_n(x) \,\nu(dx) \right)\, H_n (\xi_{k+1}) .
\end{split}
\end{equation}
Here $ f $ is a bounded measurable function.

When $ \xi $ is Gaussian, $ \{ H_n \} $ could be 
the Hermite polynomials (up to constants). 
Further, if $ \Ex[ \xi^2_k ] = t_{k} - t_{k-1} $,
then $ W $ can be a Brownian motion.

\subsection{A multi-dimensional extension}
The DIF 
for a multi-dimensional random walk
$ \mathbf{W} = (W^1,...,W^n) $, 
where 
$ W^j_{t_k} - W^j_{t_{k-1}}  = \xi^j_k $,
can be obtained 
through the expansion in 
the tensor product $ \otimes_j L^2 (\nu_j) $. 
Here the law of $ \xi^j $ is denoted by $ \nu_j $. 
Letting $ \nu_0 $ be a trivial measure, 
we get a DIF of $ (t, \bf{W}) $ as 
\begin{equation*}
\begin{split}
&f ( t_{k+1}, \mathbf{W}_{t_{k+1}} ) - f (t_k, \mathbf{W}_{t_k}) \\
&= \sum_j \frac{1}{\Ex [|\xi^j|^2 ]}
\left( \int  f (t_{k+1}, 
W^1_{t_k} + x_1,..., W^j_{t_k} + x_j ,...) x_j \,\nu_j (dx_j) \right) 
(W^j_{t_{k+1}} - W^j_{t_k})  \\
&\quad + \frac{1}{t_{k+1} -t_k}
\left( 
\int \{ f(t_{k+1}, \mathbf{W}_{t_k} +\mathbf{x} ) 
- f(t_k, \mathbf{W}_{t_k})\} \,\nu_1 \otimes \cdots \otimes \nu_n (d\mathbf{x})
\right)
\,(t_{k+1} -t_k) \\
&\qquad +\sum_{l_1 +\cdots + l_n \geq 2}
\left( \int f ( t_{k+1}, \mathbf{W}_{t_k} +\mathbf{x}) 
H^1_{l_1} (x_1) \cdots H^n_{l_n} (x_n) \,\nu_1(dx_1)\cdots
 \nu_n (dx_n) \right.\, \\
& \qquad \qquad \qquad\qquad\qquad
\left. \cdot H^1_{l_1} (\xi^1_{t_{k+1}} ) \cdots H^n_{l_n} (\xi^n_{t_{k+1}}) 
\right).
\end{split}
\end{equation*}
Here $ f $ is a bounded measurable function on $ \R^{n+1} $.

Now one sees that, in a sense, 
our DIF gives a discrete and \underline{conditioned}
chaos expansion (even when $ \xi $'s  are not Gaussian). 
The $ 0 $-th and the $ 1 $st chaoses consist of 
the main terms and the higher order terms correspond 
to the {\em correction terms}, when one wants to get the standard
It\^o formula by letting $ \Delta t := t_{k+1} - t_k $ to $ 0 $.

\subsection{DIF for solutions to stochastic difference equations}\label{SS}
Let $ \mathbf{X} $ be 
the solution of a stochastic {\em difference} equation
innovated by $ \mathbf{W} $. That is, 
\begin{equation}\label{approx_scheme}
\mathbf{X}_{t_{k+1}} = \mathbf{X}_{t_k} + 
F( \mathbf{X}_{t_k}, t_{k+1}-t_k, \mathbf{W}_{t_{k+1}} - \mathbf{W}_{t_k} )
\end{equation}
for some vector field $ F: \R^n \times \R^{n+1} \to \R^n $.
The DIF for $ (t, \mathbf{X} ) $ would be 
\begin{equation}\label{ADIF3}
\begin{split}
&f ( t_{k+1}, \mathbf{X}_{t_{k+1}} ) - f (t_k, \mathbf{X}_{t_k}) \\
&= \sum_j \frac{1}{\Ex [|\xi^j|^2 ]}
\left( \int  f (t_{k+1}, 
\mathbf{X}_{t_k} + F (\mathbf{X}_{t_k} , \Delta t,
\mathbf{x} )  ) \,x_j \,\nu_j (dx_j) \right) 
(W^j_{t_{k+1}} - W^j_{t_k})  \\
&\quad + \frac{1}{t_{k+1} -t_k}
\left( 
\int \{ f(t_{k+1}, \mathbf{X}_{t_k} + F (\mathbf{X}_{t_k} , \Delta_t,
\mathbf{x} ) )
- f(t_k, \mathbf{X}_{t_k})\} \,
\nu_1 \otimes \cdots \otimes \nu_n (d\mathbf{x})
\right)
\,(t_{k+1} -t_k) \\
&\qquad +\sum_{l_1 +\cdots + l_n \geq 2}
\left( \int f ( t_{k+1}, \mathbf{X}_{t_k} +F (\mathbf{X}_{t_k} ,\Delta t,
\mathbf{x} ) ) 
H^1_{l_1} (x_1) \cdots H^n_{l_n} (x_n) \,\nu_1(dx_1)\cdots
 \nu_n (dx_n) \right.\, \\
& \qquad \qquad \qquad\qquad\qquad
\left. \cdot H^1_{l_1} (\xi^1_{t_{k+1}} ) \cdots H^n_{l_n} (\xi^n_{t_{k+1}}) 
\right).
\end{split}
\end{equation}
When $ F (\mathbf{x},\Delta t, \mathbf{y} ) $ is affine in $ \mathbf{y} $ as
\begin{equation}\label{Maru}
F( \mathbf{x}, \Delta t, \mathbf{y}) 
= \sigma (\mathbf{x}) \mathbf{y} \sqrt{\Delta t} 
+ \mu (\mathbf{x} ) \Delta t, 
\end{equation}
where $ \sigma : \R^n \to \R^n \otimes \R^n $
and $ \mu : \R^n \to \R^n $, and when $ \mathbf{W} $ is 
a Brownian motion, 
then (\ref{approx_scheme}) can be seen as an Euler-Maruyama 
scheme of a stochastic differential equation
\begin{equation*}
d \mathbf{X} = \sigma (\mathbf{X}) \,d\mathbf{W} 
+ \mu (\mathbf{X}) \,dt.
\end{equation*} 
Still many classes, including higher order schemes 
and approximation schemes to SDE driven by L\'evy processes, 
can be also written in the form of (\ref{approx_scheme}).

\subsection{DIF for a class of weak approximation schemes}\label{WS}

For a weak approximation scheme in a Brownian cases, 
we introduce another framework.
If we define an $ n $-dimensional random walk 
$ \mathbf{W}= (W^1,...,W^n) $ by 
$ W^j_{t_k} - W^j_{t_{k-1}} = H_j (\xi_k) \sqrt{ \Delta t} $,
the (\ref{approx_scheme}) will work as a weak approximation
scheme, based on the fact that
\begin{equation*}
\sqrt{\Delta t} \left( \sum_{t_k \leq t} 
H_1 ( \xi_{t_k} ) , 
\cdots, 
\sum_{t_k \leq t} 
H_n ( \xi_{t_k} )
\right), \,\, t \geq 0
\end{equation*}
converges in law to the $ d $-dimensional Wiener process 
due to the martingale central limit theorem 
(see, e.g. \cite[Chapter 7]{EtKu}). 

In this framework, the DIF for $ (t,\bf{X}) $ becomes
\begin{equation}\label{ADIF4}
\begin{split}
&f ( t_{k+1}, \mathbf{X}_{t_{k+1}} ) - f (t_k, \mathbf{X}_{t_k}) \\
&= \sum_{j=1}^n \frac{1}{\sqrt{\Delta t}}
\left( \int  f (t_{k+1}, 
\mathbf{X}_{t_k} + F (\mathbf{X}_{t_k} , \Delta t, 
\mathbf{H}(x) ) \, ) \,H_j(x) \,\nu (dx) \right) 
(W^j_{t_{k+1}} - W^j_{t_k})  \\
&\quad + \frac{1}{\Delta t}
\left( 
\int \{ f(t_{k+1}, \mathbf{X}_{t_k} + F (\mathbf{X}_{t_k}, \Delta t,
\mathbf{H}(x) ) )
- f(t_k, \mathbf{X}_{t_k})\} \,
\nu (dx)
\right)
\,(t_{k+1} -t_k) \\
&\qquad +\sum_{l > n}
\left( \int f ( t_{k+1}, \mathbf{X}_{t_k} +F (\mathbf{X}_{t_k}, \Delta t,  
\mathbf{H} (x) ) ) 
\,H_l (x) \,\nu (dx) \right)\, H_{l} (\xi_{t_{k+1}} ),
\end{split}
\end{equation}
where we have denoted
$ \mathbf{H} (x) = ( H_1 (x),...,H_n (x) ) $.

\subsection{DIF for complete markets}\label{DIFcomplete}
If $ \sharp G =n+1 $ in section \ref{1st}, 
then $ \{ H_0,...,H_n \} $ spans the whole space and 
the correction terms in (\ref{ADIF4}) disappear. 
That is, the DIF becomes symbolically equivalent to the standard one. 
This is because the cardinality of the martingale basis
of $ \bf{W} $ is equal to the dimension of the state space
of itself, as is the case with the standard Brownian motions.

From a perspective of mathematical finance,
this property is closely related to {\em completeness} of the markets
modeled by the stochastic process.

Roughly speaking, a market is said to be 
complete if every good has a unique 
price that excludes arbitrage opportunities.
In many models in financial engineering 
the market is assumed to be complete
to avoid discussing too much 
about the utilities/preferences of individuals.
 
Discrete-time complete market models
are studied in \cite{He} \underline{using a DIF}. 

\subsection{Supplementary remarks for section \ref{DIF}}
We remark that:

\begin{enumerate}[I)]
\item This idea, namely {\em conditioned Fourier expansion of the increments} can be applied to more general cases. It can be 
``random walks on a graph/group/Polish space" , 
``discrete-time Markov chains on a manifold", or 
``general semi-martingales", etc.

\item It is also notable that our DIF holds irrespective of 
the distributions, as far as the reference measures are equivalent.

\item To the best of the author's knowledge, 
discrete It\^o's formula was pioneered by T. Szabados \cite{Sza}.

\item Discrete stochastic calculus, which have more emphasis on
chaos expansions, has been studied by many. A nice exposition 
\cite{Gzyl} on this
topic is available. 

\end{enumerate}

\section{Stochastic Numerics from the Perspective of
Discrete It\^o Calculus}\label{numerics}
\subsection{Reduction to finite difference schemes}\label{1st}

When $ \nu $ is concentrated on a finite set $ G = \{ g_1,...,g_l \} $, then 
\begin{equation*}
u^N_T (t,x) = \Ex [ f(\mathbf{X}_T) | X_t = x ], \quad t=t_0, t_1,...,t_N,
\,(t_N= T, \Delta t \equiv 1/N )
\end{equation*}
solves a \underline{finite difference} equation 
\begin{equation}\label{FDE}
\mbox{
$ ( \partial^N_t + L^N ) u = 0 $ with terminal condition $ u_T (T,x) = f(x) $,
}
\end{equation}
where 
\begin{equation*}
\partial^N_t u(t,x) = N \{ u (t, x) - u (t-1/N,x) \}
\end{equation*}
and 
\begin{equation*}
L^N u(t,\mathbf{x}) = N 
\sum_{i=1}^l  \{ u (t , \mathbf{x} + F (\mathbf{x}, 1/N,
\mathbf{H}(g_i) ) )
- u(t, \mathbf{x})\} \,
\nu (g_i).
\end{equation*}
This means that 
if $ \partial^N_t + L^N $ is {\em consistent}
\footnote{See, e.g. \cite{MR95b:65003}}
with a differential operator 
$ \partial_t + L $; i.e.
\begin{equation}\label{consistency}
\int_t^T ( \partial^N_t + L^N ) \varphi \,ds \to 
\int_t^T ( \partial_t + L ) \varphi \,ds \, \,
\mbox{as $ N \to \infty $ for any smooth $ \varphi $, }
\end{equation}
$ u^N $ converges to a smooth solution $ u $ of $  ( \partial_t + L ) u = 0 $.

Following standard arguments in the context 
of the finite difference scheme, we illustrate what is going on here.
For the solution $ v $ 
\begin{equation}\label{fds1}
(\partial^N_t + L^N ) v = \Phi, \,v(T) = 0,
\end{equation}
we have a discrete Feynman-Kac formula:
\begin{equation*}
v (t_k, x) = \sum_{l=k}^N \Ex [ \Phi (t_l, \mathbf{X}_{t_l} ) 
| \mathbf{X}_{t_k} = x ].
\end{equation*}
Suppose that the smooth solution $ u $ exists. 
Then $ u -u^N $ is a solution to (\ref{fds1}) 
with $ \Phi = \{ ( \partial^N_t + L^N ) - ( \partial_t + L ) \} u $. 
Therefore, 
\begin{equation*}
|u^N (t,x) -u(t,x) | \leq 
\sum |\Ex [\{ ( \partial^N_t + L^N ) - ( \partial_t + L ) \} u (t_l, X_{t_l}) | 
X_t = x ] |,
\end{equation*}
and the consistency (\ref{consistency}) gives the convergence.

Namely it acts as a finite difference approximation of 
the partial differential equation. 
Note that the problem of the rate of convergence in the Euler-Maruyama
scheme:
\begin{equation*}
\Ex [ f (\mathbf{X}^N_t ) ] \to \Ex [f (\mathbf{X}_T)]
\end{equation*}
reduces to the same problem in (\ref{consistency}),
which can be easily calculated in many cases.

\subsection{Completeness makes it slow}
As we remarked in section \ref{DIFcomplete}, the cases 
where $ \sharp G = 1 + $ the dimension of the state space
is of special interest since it serves as a complete market model.
However, if one considers them to be a discretization,
by the Euler-Maruyama scheme,
of a continuous-time (complete market) model, one 
is obliged to pay some costs. 
\begin{thm}[\cite{He}]\label{zero}
The unique prices of European claims in discrete-time 
complete markets 
converge
to the ones in the corresponding continuous complete market 
as the time-step $ \Delta t $ tends to $ 0 $. 
The order of convergence is at least $ \sqrt{\Delta t} $ 
and {\bf it cannot be improved} when 
$ n \geq 2 $.
\end{thm}

Roughly speaking, this is because the set $ G $ with $ \sharp G = n+1 $
{\bf cannot support any} $ n $-dimensional random variable
which has the same moments up to degree three 
with the increment of 
$ n $-dimensional Brownian motion, when $ n \geq 2 $.
For details, see \cite{He}.

\subsection{``Infinite" difference scheme}\label{high}
The argument in section \ref{1st} is still valid for
a general $ \nu $ by putting 
\begin{equation}\label{semi2}
L^N f (\mathbf{x}) = N \int f 
(\mathbf{x} + F (\mathbf{x} , N^{-1},
\mathbf{H}(y) ) ) \,\nu (dy)
\end{equation}
for the SDE's in section \ref{WS}.

We propose the following implementations.
\begin{enumerate}
\item Let $ \nu $ be the Lebesgue measure on $ [0,1] $.
\item Let $ F $ be as (\ref{Maru}) (Euler-Maruyama).
\item Take the Walsh system as an ONB. 
Here by the Walsh system we mean 
the group generated by 
\begin{equation*}
\tau_{n} (x) = 
\begin{cases}
+1 & \mbox{ if $ 2k \leq 2^n x < 2k+1 $ for 
$ k = 0,1,...,2^{n-1} -1 $,} \\
-1 & \mbox{otherwise},
\end{cases}
\end{equation*}
$ n=1,2, \cdots$.
Note that $  \tau_1,...,\tau_n,... $ are nothing but a Bernoulli sequence.
We can take, for example, 
\begin{equation*}
\mathbf{H} = (\tau_1, \tau_2, \tau_3, \tau_1\tau_2\tau_3, 
\tau_4, \tau_5, \tau_1 \tau_2 \tau_4, \tau_1\tau_2 \tau_5, 
\tau_1 \tau_2 \tau_3 \tau_4 \tau_5, .... ).
\end{equation*}
(Just avoid using those with the even-number length.)

\item Simulate the path $\mathbf{X} $ by a 
Monte-Carlo/Quasi Monte-Carlo uniform sequence in $ [0,1)^N $.
\end{enumerate}
This method is meant to be a Monte-Carlo simulation scheme 
of high-dimensional finite difference scheme. 
It is almost dimension-free. 
In fact, the dimension $ n $ can be effectively 
very large; around 3000, as is reported in \cite{Sawai}. 

\subsection{Supplementary remarks for section \ref{numerics}}
The correspondence of the general Markov chain approximation by
(\ref{approx_scheme}) to finite-difference schemes can be generalized
to non-linear cases. This generalization includes
the Kushner's correspondence (see e.g. \cite{Kushner}).

\end{document}